\documentclass[10pt]{article}

\usepackage{amsmath}
\usepackage{amssymb}
\usepackage{latexsym}
\usepackage{amsthm}

\newcommand{\init}{\big\vert_{t = 0}}
\newcommand{\abs}[1]{\left\vert #1 \right\vert}
\newcommand{\bigabs}[1]{\bigl\vert #1 \bigr\vert}

\newcommand{\norm}[1]{\left\Vert #1 \right\Vert}
\newcommand{\bignorm}[1]{\bigl\Vert #1 \bigr\Vert}

\newcommand{\spacetimenorm}[3]{\norm{#1}_{H^{#2,#3}}}

\newcommand{\Spacetimenorm}[3]{\norm{#1}_{\scrH^{#2,#3}}}

\newcommand{\Sobnorm}[2]{\norm{#1}_{H^{#2}}}
\newcommand{\bigSobnorm}[2]{\bignorm{#1}_{H^{#2}}}
\newcommand{\datanorm}[3]{\norm{(#1,#2)}_{(#3)}}

\newcommand{\Lpnorm}[2]{\norm{#1}_{L^{#2}}}

\newcommand{\twonorm}[2]{\norm{#1}_{L^2#2}}

\newcommand{\inftynorm}[2]{\norm{#1}_{L^\infty#2}}
\newcommand{\biginftynorm}[2]{\bignorm{#1}_{L^\infty#2}}

\newcommand{\Bdualnorm}[1]{\norm{#1}_{L_{\xi}^{2}(L_{\tau}^{\infty})}}
\newcommand{\bigBdualnorm}[1]{\bignorm{#1}_{L_{\xi}^{2}(L_{\tau}^{\infty})}}

\newcommand{\C}{\mathbb{C}}
\newcommand{\D}{\mathcal{D}}
\newcommand{\half}{\frac{1}{2}}
\newcommand{\scrH}{\mathcal{H}}

\newcommand{\R}{\mathbb{R}}
\newcommand{\X}{\mathcal{X}}
\newcommand{\Y}{\mathcal{Y}}

\newcommand{\Schwartz}{\mathcal{S}}
\newcommand{\Fourier}{\mathcal{F}}
\newcommand{\Distribution}{\mathcal{D}'}

\newcommand{\hypwt}[2]{\bigabs{ \abs{#1} - \abs{#2} }}

\DeclareMathOperator{\supp}{supp}

\newtheorem{theorem}{Theorem}

\newtheorem{lemma}{Lemma}

\theoremstyle{definition}

\theoremstyle{remark}

\newtheorem*{plainremark}{Remark}
\newtheorem*{plainremarks}{Remarks}
\title{On an estimate for the wave equation and applications to nonlinear problems}
\author{Sigmund Selberg\\ Department of Mathematics\\ Johns Hopkins University
\\Baltimore, MD 21218}
\date{}
\begin{document}
\maketitle
{
\renewcommand{\thefootnote}{}
}
\begin{abstract}
We prove estimates for solutions of the Cauchy problem
for the inhomogeneous wave equation on $\R^{1+n}$ in a class of
Banach spaces whose norms only depend on the size of the space-time Fourier
transform. The estimates are local in time, and this
allows one, essentially, to replace
the symbol of the wave operator, which vanishes on the light cone in Fourier space,
with an inhomogeneous symbol, which can be inverted.
Our result improves earlier estimates of this type
proved by Klainerman-Machedon \cite{Kl-Ma1,Kl-Ma3}.
As a corollary, one obtains a rather general result concerning
local well-posedness of nonlinear wave equations, which
was used extensively in the recent article \cite{Kl-Se}.
\end{abstract}
\section{Introduction}
Consider the Cauchy problem for the wave equation on $\R^{1+n}$,
\begin{equation}\label{CauchyProblem}
  \square u = F, \quad u \init = f, \quad \partial_t u \init = g,
\end{equation}
where $\square = -\partial_t^2 + \Delta$ is the wave operator.

The purpose of this
note is to prove estimates for $u$ in a certain class of Banach spaces whose
norms only depend on the size of the space-time Fourier transform. The estimates
are local with respect to time $t$, and this allows one, essentially, to replace
the symbol of the wave operator, which vanishes on the light cone in Fourier space,
with an inhomogeneous symbol, which can be inverted. This idea originates in the
work of Bourgain \cite{B} on the Schr\" odinger and KdV equations, and was later
simplified by Kenig-Ponce-Vega \cite{KPV} in their work on KdV.
Following this, Klainerman-Machedon \cite[Lemma 4.3]{Kl-Ma1},
\cite[Lemma 1.3]{Kl-Ma3} proved estimates of this type for the wave equation;
see also Klainerman-Tataru \cite{Kl-Ta}.

The improvement in our result compared to \cite{KPV,Kl-Ma1,Kl-Ma3}
lies mainly in showing, as suggested in \cite[Remark 1.8]{Kl-Ta}, that for sufficiently
small $\varepsilon > 0$, the norm of the inhomogeneous part of the solution,
restricted to the time slab $[0,T] \times \R^n$, is $O(T^\varepsilon)$ as $T \to 0$,
at the expense of a loss of essentially $\varepsilon$ derivatives.
As shown in the recent article \cite{Kl-Se}, this allows one
to remove the assumption of small-norm data in the well-posedness results
proved in \cite{Kl-Ma1,Kl-Ma3,Kl-Ma5,Kl-SeWM,Kl-Ta}.

As an application of our estimate, we also prove a simple but useful
result concerning local well-posedness of nonlinear wave equations,
which is used extensively in \cite{Kl-Se}.
\section{The main estimate}\label{Main}
We are interested in finding complete subspaces $\X^s$ of
\begin{equation}\label{GlobalEnergySpace}
  C_b(\R,H^s) \cap C^1_b(\R,H^{s-1})
\end{equation}
such that solutions of \eqref{CauchyProblem} satisfy estimates of the type
\begin{equation}\label{BasicEstimate}
  \norm{u}_{\X^s_T} \le C \datanorm{f}{g}{s} + C_{T,\varepsilon}
  \norm{\widetilde \square ^{-1} F^\varepsilon}_{\X^s}
\end{equation}
for all $0 < T < 1$ and $\varepsilon \ge 0$.
We also want
\begin{equation}\label{ConstantLimit}
  \lim_{T \to 0} C_{T,\varepsilon} = 0
\end{equation}
when $\varepsilon$ is strictly positive.

Precise definitions will be supplied presently. For the moment suffice it to
say that $\X^s_T$ stands for the restriction to the time slab $[0,T] \times \R^n$,
$\widetilde \square$ may be thought of as an inhomogeneous and invertible
version of the wave operator $\square$,
$F^\varepsilon$ is $F$ with a certain operator of order $\varepsilon$ applied to it, and
$$
  \datanorm{f}{g}{s} = \Sobnorm{f}{s} + \Sobnorm{g}{s-1}
$$
with $H^s$ the usual Sobolev space. 

We use coordinates $(t,x)$ on $\R^{1+n}$.
The Fourier transform of $f(x)$ [resp. $u(t,x)$] is denoted $\widehat f(\xi) = \Fourier f(\xi)$
[resp. $\widehat u(\tau,\xi) = \Fourier u(\tau,\xi)$].
For any $\alpha \in \R$ we define pseudodifferential operators
$\Lambda^{\alpha}$, $\Lambda_{+}^{\alpha}$ and $\Lambda_{-}^{\alpha}$ by
\begin{align*}
  \widehat{\Lambda^{\alpha} f}(\xi) &= \bigl( 1 + \abs{\xi}^{2}
  \bigr)^{\alpha/2} \widehat{f}(\xi),
  \\
  \widehat{\Lambda_{+}^{\alpha} u}(\tau,\xi) &= \bigl( 1 + \tau^{2} + \abs{\xi}^2
  \bigr)^{\alpha/2} \widehat{u}(\tau,\xi),
  \\
  \widehat{\Lambda_{-}^{\alpha} u}(\tau,\xi) &=
  \left( 1 + \frac{(\tau^2 - \abs{\xi}^2)^{2} }
  {1 + \tau^{2} + \abs{\xi}^2 } \right)^{\alpha/2}
   \widehat{u}(\tau,\xi),
\end{align*}
Observe that the Fourier symbol of $\Lambda_-^\alpha$ is comparable
to $\bigl(1 + \hypwt{\tau}{\xi} \bigr)^\alpha$. The operator $\widetilde \square$
in \eqref{BasicEstimate} is just $\Lambda_+ \Lambda_-$, and
$F^\varepsilon = \Lambda_-^\varepsilon F$.

The Sobolev and Wave Sobolev spaces $H^s$ and $H^{s,\theta}$ are given
by the norms
\begin{align*}
  \Sobnorm{f}{s} &= \twonorm{\Lambda^s f}{(\R^n)},
  \\
  \spacetimenorm{u}{s}{\theta} &= \twonorm{\Lambda^s \Lambda_-^\theta u}{(\R^{1+n})}.
\end{align*}
For the basic properties of the latter, see, e.g., \cite{Kl-Se}.
In particular, we shall use the fact that $H^{s,\theta}$ embeds in $C_b(\R,H^s)$
when $\theta > \half$.
Associated to $H^{s,\theta}$ is the space $\scrH^{s,\theta}$ with norm
$$
  \Spacetimenorm{u}{s}{\theta} = \spacetimenorm{u}{s}{\theta}
  + \spacetimenorm{\partial_t u}{s-1}{\theta}
  \sim \twonorm{\Lambda^{s-1} \Lambda_+ \Lambda_-^\theta u}{}.
$$
By the above, $\scrH^{s,\theta}$ embeds in \eqref{GlobalEnergySpace}
for $\theta > \half$.

In general, if a Banach space $\X^s$ embeds in \eqref{GlobalEnergySpace}
then it makes sense to restrict its elements to any interval $I \subset \R$.
The resulting restriction space is denoted $\X^s_I$. It is always possible to define
a norm on this space which makes it complete. Indeed, $\X^s_I$
is the quotient space $\X^s / \!\!\sim_{I}$, where $\sim_I$ is the equivalence
relation
$$
  u \sim_{I} v \iff \text{$u(t) = v(t)$ for all $t \in I$}.
$$
Since $\X^{s}$ embeds in \eqref{GlobalEnergySpace}, the equivalence
classes are closed sets in $\X^s$, so the quotient space
is complete when equipped with the norm
$$
  \norm{u}_{\X^s_I} = \inf_{v \sim_{I} u} \norm{v}_{\X^s}.
$$
If $I = [0,T]$, we always write $\X^s_T$ instead of $\X^s_I$.

We now state the precise result.
\begin{theorem}\label{BasicConditionsTheorem} Let $\X^{s}$ be a Banach
space such that
\begin{enumerate}
  \item $\X^s$ embeds in $\scrH^{s,\theta}$ for some $\theta > \half$,
  \item  $\abs{\widehat u} \le \abs{\widehat v}
  \implies \norm{u}_{\X^s} \le \norm{v}_{\X^s}$,
  \item there exists $\gamma < 2$ such that\footnote{We use the notation
  $\norm{v(\tau,\xi)}_{L_{\xi}^{2}(L_{\tau}^{\infty})}
  = \bigl( \int \inftynorm{v(\cdot,\xi)}{}^2 \, d\xi \bigr)^{1/2}$.}
  $$\norm{u}_{\X^{s}} \lesssim \bignorm{ \Fourier \Lambda^{s-1} \Lambda_{+}
  \Lambda_{-}^{\gamma} u (\tau,\xi) }_{L_{\xi}^{2}(L_{\tau}^{\infty})}$$
  for all $u$.
\end{enumerate}
Let $\varepsilon \ge 0$. Then for all $(f,g) \in H^s \times H^{s-1}$ and
$F \in \widetilde \square \Lambda_{-}^{-\varepsilon} ( \X^{s} )$,
there is a unique $u \in C(\R,H^s) \cap C^1(\R,H^{s-1})$ which solves \eqref{CauchyProblem},
and the estimate \eqref{BasicEstimate} holds for all $0 < T < 1$.
Moreover, if $\varepsilon > 0$, then \eqref{ConstantLimit} holds.
\end{theorem}
\begin{plainremarks}
(1) Estimate \eqref{BasicEstimate} in the special case $\X^s = \scrH^{s,\theta}$
and $\varepsilon = 0$ was proved in \cite{Kl-Ma1}.

\medskip

\noindent
(2) The proof gives something stronger: for all $T > 0$, there
is a linear operator $W_T$ such that
$$
  \square W_T F = F \quad \text{on} \quad [0,T] \times \R^n,
$$
with vanishing initial data at $t = 0$,
and $W_T$ is bounded from
\begin{equation}\label{Yspace}
  \Y^{s,\varepsilon} =
  \widetilde \square \Lambda_-^{-\varepsilon} (\X^s)
\end{equation}
into $\X^s$
for all $\varepsilon \ge 0$. Moreover, when $\varepsilon > 0$,
the operator norm
$$
  \norm{W_T}_{\Y^{s,\varepsilon} \to \X^s} \to 0 \quad \text{as}
  \quad T \to 0.
$$

\medskip
\noindent
(3) The reason for the condition $\gamma < 2$ is as follows.
The proof of the theorem shows that
$$
  C_{T,\varepsilon} \le C T^{\min\{\alpha \varepsilon, \delta\}},
$$
where
$$
  \delta = 2 - \gamma - \alpha/2
  + \min \bigl\{ 0, \alpha(\theta + \varepsilon - 1) \bigr\}
$$
and $0 \le \alpha \le 1$ can be chosen at will. The constant
$C$ is independent of $\varepsilon$ and $\alpha$.
In the typical applications
$\gamma$ is close to 1; see \cite{Kl-Se}.

\medskip
\noindent
(4) $s$ plays no role. Indeed, if $\X^s$ satisfies the hypotheses of the theorem,
then so does $\X^{s'} = \Lambda^{s-s'} \X^s$ for all $s'$.
\end{plainremarks}
\section{Applications}
For $t \in \R$ we denote by $\tau_t$ the time-translation
operator $\tau_t u = u(\cdot + t,\cdot)$, and for any interval
$I \subseteq \R$ we denote restriction to the time-slab
$I \times \R^n$ by $\vert_I$.
 
Suppose $\X^s$ satisfies the hypotheses of
Theorem \ref{BasicConditionsTheorem}, is
invariant under time-translation\footnote{That is to say, $\tau_t$ is an
isomorphism of $\X^s$ for all $t$.},
and for all $\phi \in C_{c}^{\infty}(\R)$, the multiplication map
$u \mapsto \phi(t) u(t,x)$ is bounded from $\X^s$ into itself.

Consider a system of wave equations on $\R^{1+n}$ of the form
\begin{equation}\label{GenericSystem}
  \square u = \mathcal N(u),
\end{equation}
where $u$ takes values in $\R^{N}$ and $\mathcal N : \X^s \to \Distribution$
is (i) time-translation invariant ($\mathcal N$ commutes with $\tau_t$); (ii)
local in time\footnote{By this we mean that if $u \vert_I = v \vert_I$,
where $I$ is an open interval, then $\mathcal N(u) \vert_I
= \mathcal N(v) \vert_I$.}; and (iii) satisfies $\mathcal N(0) = 0$.

Furthermore, we assume that for some $\varepsilon > 0$,
$\mathcal N$ satisfies, with notation as in \eqref{Yspace},
a Lipschitz condition
\begin{equation}\label{EpsilonGenericDifferenceEstimate}
  \norm{\mathcal N (u) - \mathcal N (v)}_{\Y^{s,\varepsilon}}
  \le A\bigl( \max\{\norm{u}_{\X^s},\norm{v}_{\X^s} \} \bigr) \norm{u-v}_{\X^s}
\end{equation}
for all $u,v \in \X^{s}$, where $A$ is a continuous function.
\begin{theorem}\label{SpecialWPTheorem}
Under the above assumptions, the Cauchy problem
for \eqref{GenericSystem} is locally well-posed for initial data in
$H^{s} \times H^{s-1}$, in the following sense:
\begin{itemize}
 \item[(I)] {\bfseries(Local existence)}
  For all $(f,g) \in H^{s} \times H^{s-1}$ there exist a $T > 0$
  and a $u \in \X^s_{T}$ which solves \eqref{GenericSystem} on
  $S_{T} = (0,T) \times \R^{n}$ with initial data $(f,g)$.
  Moreover, $T$ can be chosen to depend continuously on 
  $\Sobnorm{f}{s} + \Sobnorm{g}{s-1}$.
  \item[(II)] {\bfseries(Uniqueness)} If $T > 0$ and $u,u' \in \X^s_T$
  are two solutions of \eqref{GenericSystem} on $S_{T}$
  with the same initial data $(f,g)$, then $u = u'$.
  \item[(III)] {\bfseries(Continuous dependence on initial data)}
  If, for some $T > 0$, $u \in \X^{s}_{T}$ solves \eqref{GenericSystem} on $S_{T}$
  with initial data $(f,g)$, then for all $(f',g') \in H^{s} \times H^{s-1}$
  sufficiently close to $(f,g)$, there exists a $u' \in \X^{s}_{T}$
  which solves \eqref{GenericSystem} on $S_{T}$ with initial data $(f',g')$,
  and
  $$
    \norm{u - u'}_{\X^{s}_{T}} \le C \datanorm{f-f'}{g-g'}{s}.
  $$
\end{itemize}
If, moreover, $\mathcal N$ is $C^\infty$ as a map from $\X^s$ into
$\Y^{s,\varepsilon}$, then:
\begin{itemize}
\item[(IV)] {\bfseries(Smooth dependence on initial data)}
Suppose $u \in \X^s_T$ solves \eqref{GenericSystem} on
$S_T$ for some
$T > 0$, with initial data $(f,g) \in H^s \times H^{s-1}$, and that
$(f_\delta,g_\delta)$ is a smooth perturbation of the initial data,
i.e.,
$$
  \delta \mapsto (f_\delta,g_\delta),
  \qquad \R \to H^s \times H^{s-1}
$$
is $C^\infty$ and takes the value $(f,g)$ at $\delta = 0$. Let $u_\delta$
be the corresponding solution of \eqref{GenericSystem} (by (III),
$u_\delta \in \X^s_T$ for $\abs{\delta} < \delta_0$).
Then the map $\delta \mapsto u_\delta$ from $[-\delta_0,\delta_0]$ into
$\X^s_T$ is $C^\infty$.
\end{itemize}
We write $\X^\sigma = \Lambda^{s-\sigma} \X^s$
and $\Y^{\sigma,\varepsilon} = \Lambda^{s-\sigma}\Y^{s,\varepsilon}$
for $\sigma \in \R$.
If for all $\sigma > s$ there is a continuous
function $A_\sigma$ such that
\begin{equation}\label{HigherOrderEstimate}
  \norm{\mathcal N (u) }_{\Y^{\sigma,\varepsilon}}
  \le A_\sigma \bigl( \norm{u}_{\X^s} \bigr) \norm{u}_{\X^\sigma}
\end{equation}
for all $u \in \X^s \cap \X^{\sigma}$, then:
\begin{itemize}
  \item[(V)]
  {\bfseries(Persistence of higher regularity)}
  If $\sigma > s$ and $u \in \X^s_T$ solves \eqref{GenericSystem} on
  $S_T$ with initial data $(f,g) \in H^\sigma \times H^{\sigma-1}$
  for some $T > 0$, then
  $$
    u \in C\bigl([0,T],H^{\sigma}\bigr) \cap 
    C^{1}\bigl([0,T],H^{\sigma-1}\bigr).
  $$
\end{itemize}
\end{theorem}
\begin{plainremark}
Typically, proving that $\mathcal N$ is $C^\infty$ is no
harder than proving it is locally Lipschitz. As an example,
relevant for wave maps, consider
$$
  \mathcal N(u) = \Gamma(u) Q_0(u,u),
$$
where $u$ is real-valued, $\Gamma : \R \to \R$ is $C^\infty$ and
$Q_0$ is the bilinear ``null form''
$$
  Q_0(u,v) = -\partial_t u \, \partial_t v + \nabla_x u \cdot \nabla_x v.
$$
Fix $s > \frac{n}{2}$ and set
$\X = \scrH^{s,\theta}$, $\X' = H^{s,\theta}$ and
$\Y = \widetilde \square \Lambda_-^{-\varepsilon}(\X)
= H^{s-1,\theta+\varepsilon-1}$. One can show,
for appropriate $\theta > \half$ and $\varepsilon > 0$, that
\begin{itemize}
  \item[(i)]
  $Q_0$ is bounded, hence $C^\infty$, from $\X \times \X$ to $\Y$;
  \item[(ii)]
  $\Phi_\Gamma : u \mapsto \Gamma(u)-\Gamma(0)$ is a locally bounded map of $\X'$ to itself;
  \item[(iii)]
  multiplication is bounded, hence $C^\infty$, from $\X' \times \Y$
  to $\Y$.
\end{itemize}
It remains to prove that $\Phi_\Gamma$ is $C^\infty$ as a map of
$\X'$, but this follows from (ii) (which is valid for any $C^\infty$
function $\Gamma$) and the fact that $\X'$ is an algebra.
Indeed, since
$$
  \Gamma(v) - \Gamma(u) =
  \int_0^1 \Gamma' \bigl(u + t(v-u)\bigr)(v-u)
  \, dt,
$$
$\Phi_\Gamma$ is locally Lipschitz. Then, since
$$
  \Gamma(v) - \Gamma(u) - \Gamma'(u)(v-u)
  = \int_0^1 \left\{ \Gamma' \bigl(u + t(v-u)\bigr) - \Gamma'(u) \right\}
  \, dt \cdot (v-u)
$$
we have
$$
  \norm{\Gamma(v) - \Gamma(u) - \Gamma'(u)(v-u)}_{\X'}
  \le C(u,v) \norm{v-u}_{\X'},
$$
where
$$
  C(u,v) \le C \int_0^1
  \norm{\Gamma' \bigl(u + t(v-u)\bigr) - \Gamma'(u)}_{\X'}
  \, dt.
$$
But by the above, $\Phi_{\Gamma'}$ is locally Lipschitz,
so $C(u,v) = O(\norm{v-u}_{\X'})$.
Thus $\Phi_\Gamma$ is $C^1$, and by induction $C^\infty$.
\end{plainremark}
\section{Abstract local well-posedness}\label{AbstractLWP}
Theorem \ref{SpecialWPTheorem} is conveniently
proved in an abstract setting, which we discuss here.

Observe that if $\X^s$ is a space satisfying the assumptions of
the last section, and $\mathcal I$ denotes the set of compact intervals
$I \subseteq \R$, then the family $\{\X^s_I\}_{I \in \mathcal I}$
has the following properties:
\begin{itemize}
  \item[(S1)]
  $\X^s_{I}$ embeds in $C(I,H^{s}) \cap C^{1}(I,H^{s-1})$ for all $I 
  \in \mathcal I$;
  \item[(S2)]
  the solution of $\square u = 0$ with
  $(u,\partial_{t} u) \init = (f,g) \in H^s \times H^{s-1}$ satisfies
  $$
    \norm{u}_{\X^{s}_{T}} \le C \datanorm{f}{g}{s}
  $$
  for all $0 < T < 1$ and all $(f,g)$;
  \item[(S3)]
  $\tau_t$ is an isometry from $\X^{s}_{I}$
  onto $\X^{s}_{-t + I}$ for all $t \in \R$ and $I \in \mathcal I$;
  \item[(S4)]
  if $I \subseteq J$, then $\vert_{I} : \X^{s}_{J} \to \X^{s}_{I}$
  is norm decreasing;
  \item[(S5)] whenever $I$ and $J$
  are two overlapping intervals ($I \cap J$ has nonempty 
  interior) and $u \in \X^{s}_{I}$, $v \in \X^{s}_{J}$ agree on the overlap
  ($u(t) = v(t)$ for all $t \in I \cap J$), then
  $$
    \norm{w}_{\X^{s}_{I \cup J}} \le C_{I,J} \left( 
    \norm{u}_{\X^{s}_{I}} + \norm{v}_{X^{s}_{J}} \right),
  $$
  where $w(t) = u(t)$ for $t \in I$ and $w(t) = v(t)$ for $t \in J$.
\end{itemize}
Note that (S2) holds by Theorem \ref{BasicConditionsTheorem};
we write $\X^{s}_{T}$ instead of $\X^{s}_{[0,T]}$.
(S5) holds by the assumption that 
$u \mapsto \phi(t) u(t,x)$ is bounded on $\X^s$
for all $\phi \in C_{c}^{\infty}(\R)$.

We now consider an abstract family of spaces with these properties.
\begin{theorem}\label{AbstractWPTheorem}
Let $s \in \R$. Let $\bigl\{ \X^{s}_{I} \bigr\}_{I \in \mathcal I}$
be a family of Banach spaces satisfying (S1--5).
Consider the system \eqref{GenericSystem}, where $\mathcal N$ 
is an operator which
\begin{itemize}
  \item[(N1)]
  maps $\X^{s}_{[a,b]}$ into $\Distribution\bigl( (a,b) 
  \times \R^{n} \bigr)$ for all $-\infty < a < b < \infty$;
  \item[(N2)]
  is time-translation invariant: $\mathcal N \circ \tau_t = \tau_t \circ \mathcal N$;
  \item[(N3)]
  is local in time: $\mathcal N(u\vert_{[a,b]}) = \mathcal N(u) \vert_{(a,b)}$
  whenever $[a,b] \subset I$ and $u \in \X^s_I$;
  \item[(N4)]
  satisfies $\mathcal N(0) = 0$.
\end{itemize}
Assume further that for all $0 < T < 1$ and $u \in 
\X^{s}_{T}$, there exists $v \in \X^{s}_{T}$ (necessarily unique) which solves
$\square v = \mathcal N(u)$ on $S_T = (0,T) \times \R^{n}$ with vanishing 
initial data at $t = 0$; we write $v = W \mathcal N(u)$ 
and assume that
\begin{equation}\label{NonlinearEstimate}
  \norm{W \mathcal N(u)}_{\X^{s}_{T}}
  \le C_{T} A \bigl( \norm{u}_{\X^{s}_{T}} \bigr) \norm{u}_{\X^{s}_{T}}
\end{equation}
and, more generally,
\begin{equation}\label{NonlinearDifferenceEstimate}
  \norm{W \bigl( \mathcal N(u) - \mathcal N(v) \bigr)}_{\X^{s}_{T}}
  \le C_{T} A \bigl( \max \{ \norm{u}_{\X^{s}_{T}}, 
  \norm{v}_{\X^{s}_{T}} \} \bigr) \norm{u - v}_{\X^{s}_{T}}
\end{equation}
for all $0 < T < 1$ and $u, v \in \X^{s}_{T}$, where
\begin{equation}\label{CTDecay}
  \lim_{T \to 0^{+}} C_{T} = 0
\end{equation}
and $A$ is a continuous function.

Then the system \eqref{GenericSystem} is
locally well-posed for initial data in $H^{s} \times H^{s-1}$,
in the sense that properties
(I--III) of Theorem \ref{SpecialWPTheorem} hold.
\end{theorem}
\begin{plainremarks}
(1) (\emph{Lifespan.})
By local existence (I), and properties (S4) and (N3), for given
$(f,g)$ the set $E(f,g)$ consisting of all $T>0$ for which
there exists $u \in \X^{s}_{T}$
solving \eqref{GenericSystem} on $S_{T}$ with data $(f,g)$,
is a nonempty interval. Then, by local existence (I) and uniqueness (II),
as well as properties (S3--5) and (N2,3), it follows that this interval
is open. Moreover, by continuous dependence on initial data (III), the 
\emph{lifespan} $T^{*} = \sup E(f,g)$ is a lower
semicontinuous function of $(f,g)$.

\medskip
(2) (\emph{Higher regularity.})
Set $\X^\sigma_T = \Lambda^{s-\sigma} \X^s_T$.
If, for any $\sigma > s$, there is a continuous $A_\sigma$ such that
\begin{equation}\label{HigherOrderNonlinearEstimate}
  \norm{W \mathcal N(u)}_{\X^{\sigma}_{T}}
  \le C_{T} A_\sigma \bigl( \norm{u}_{\X^{s}_{T}} \bigr)
  \norm{u}_{\X^{\sigma}_{T}}
\end{equation}
and
\begin{multline}\label{HigherOrderNonlinearDifferenceEstimate}
  \norm{W \bigl( \mathcal N(u) - \mathcal N(v)
  \bigr)}_{\X^{\sigma}_{T}}
  \le C_{T} A_\sigma \bigl( \max \{ \norm{u}_{\X^{s}_{T}}, 
  \norm{v}_{\X^{s}_{T}} \} \bigr) \norm{u - v}_{\X^{\sigma}_{T}}
  \\
  + C_{T} A_\sigma \bigl( \max \{ \norm{u}_{\X^{\sigma}_{T}}, 
  \norm{v}_{\X^{\sigma}_{T}} \} \bigr) \norm{u - v}_{\X^{s}_{T}}
\end{multline}
for all $u,v \in \X^s_T \cap \X^\sigma_T$, where $C_T$ is the constant
appearing in \eqref{NonlinearEstimate} and
\eqref{NonlinearDifferenceEstimate}, then (V) of Theorem
\ref{SpecialWPTheorem} holds.

\medskip
(3) (\emph{Smooth dependence on data.})
Suppose $S = W \mathcal N$ is not just Lipschitz but $C^\infty$
as a map of $\X^s_T$.
Recall that the $k$-th derivative $S^{(k)}(u)$, $u \in \X^s_T$,
is a $k$-linear map from $\X^s_T \times \cdots \times \X^s_T$
into $\X^s_T$; let $\norm{S^{(k)}(u)}_{(T)}$ denote
its operator norm. In view of \eqref{HigherOrderNonlinearDifferenceEstimate},
\begin{equation}\label{FirstDerBd}
  \bignorm{S'(u)}_{(T)} \le C_T A\bigl( \norm{u}_{\X^s_T} \bigr).
\end{equation}
Suppose there exist, for $k = 2,3,\dots$, $B_k$ continuous and increasing so that
\begin{equation}\label{HigherDerBd}
  \sup_{0 < T < 1} \sup_{\norm{u}_{\X^s_T} \le R}
  \bignorm{S^{(k)}(u)}_{(T)} \le B_k(R).
\end{equation}
Then (IV) of Theorem \ref{SpecialWPTheorem} holds.

\medskip
\noindent
(4) If $\mathcal N$ is multilinear, then so is $W \mathcal N$,
so if the latter is bounded on some Banach space, it is trivially
$C^\infty$ on that space. Thus, by the previous remark,
the dependence on initial data is $C^\infty$.

In this connection, we mention an
interesting observation due to Keel-Tao \cite[Section 8]{Ke-Ta},
concerning the feasibility of proving well-posedness for wave maps
in the critical data space by an iteration argument. For simplicity
we take $n = 2$, but this is not essential.

A wave map $u : \R^{1+2} \to S^1 \subseteq \C$ satisfies the equation
\begin{subequations}\label{WMCauchy}
\begin{equation}\label{SphereWM}
  \square u + u (\partial_\mu u \cdot \partial^\mu u) = 0,
\end{equation}
where $\cdot$ is the Euclidean inner product on $\R^2 = \C$.
Consider initial data
\begin{equation}\label{WMdata}
  u \init = 1, \qquad \partial_t u \init = ig,
\end{equation}
\end{subequations}
where $i$ is the imaginary unit and $g \in L^2$ is real-valued.

Observe that if \eqref{WMCauchy} is well-posed for
$g \in L^2$, then the solutions stay on $S^1$ and hence
are wave maps. This is certainly true for any smooth solution,\footnote{If
$u : \R^{1+2} \to \R^2$ is a smooth solution of \eqref{WMCauchy}
then $\phi = u \cdot u - 1$ solves the linear equation
$\square \phi + \phi(\partial_\mu u \cdot \partial^\mu u) = 0$ with vanishing
initial data, so by uniqueness, $\phi$ must vanish.}
and therefore true in general by an approximation argument, using the
continuous dependence on initial data.

If \eqref{WMCauchy} could be proved well-posed for
$g \in L^2$ by an iteration argument in
some Banach space, the dependence on the initial data
would necessarily be $C^\infty$, since in this case the operator
$\mathcal N$ is trilinear.

As it turns out, \eqref{WMCauchy} is well-posed for $g \in L^2$,
but the dependence on the data is not even $C^2$.
In fact, since $\theta \mapsto e^{i\theta}$, $\R \to S^1$ is a geodesic,
the solution of \eqref{WMCauchy} is given by
$u = e^{i v}$, where $\square v = 0$ with initial data $(0,g)$; clearly
$(u,\partial_t u)$ belongs to $C(\R,\dot H^1 \times L^2)$ and depends
continuously on $g$.
 
Thus $u_\varepsilon = e^{i \varepsilon v}$ solves
\eqref{WMCauchy} with $g$ replaced by $\varepsilon g$.
But since $\dot H^1$ is not an algebra, one would not expect the map
$\varepsilon \to u_\varepsilon \vert_{t = 1}$, $\R \to \dot H^1$ to be
twice differentiable at $\varepsilon = 0$ for all choices of $g$,
and indeed it is not, as proved in \cite[Proposition 8.3]{Ke-Ta}.
\end{plainremarks}
\section{Proof of Theorem \ref{BasicConditionsTheorem}}
\label{ProofOfBasicConditionsTheorem}
In this section, (i--iii) refer to the hypotheses
of the theorem. With notation as in \eqref{Yspace}, observe that
by (i),
\begin{equation}\label{Yembedding}
  \Y^{s,\varepsilon} \subseteq H^{s-1,\theta+\varepsilon-1}
\end{equation}
for all $\varepsilon \ge 0$, and that $\Y^{s,\varepsilon} \subseteq \Y^{s,0}$.

We shall use the fact that the symbol of $\Lambda_-$ is comparable
to $1 + \hypwt{\tau}{\xi}$. More precisely, there is a constant $c > 0$ such that
\begin{equation}\label{LambdaMinusSymbolEquivalence}
  c^{-1} \bigl( 1 + \hypwt{\tau}{\xi} \bigr)
  \le
  \left( 1 + \frac{ (\tau^2 - \abs{\xi}^2)^{2} }
  {1 + \abs{\tau}^{2} + \abs{\xi}^2 } \right)^{\half}
  \le
  c \bigl( 1 + \hypwt{\tau}{\xi} \bigr).
\end{equation}

Fix two bump functions $\chi, \phi \in C_{c}^{\infty}(\R)$ such that
$0 \le \chi, \phi \le 1$,
$\chi = 1$ on $[-2,2]$, and
$\phi = 1$ on $[-2c,2c]$ with support in $[-4c,4c]$,
where $c$ is the constant in \eqref{LambdaMinusSymbolEquivalence}.

Fix $0 \le \alpha \le 1$. Given $0 < T < 1$, write, for any $F \in H^{s-1,\theta-1}$,
$$
  F = \phi\bigl( T^\alpha \Lambda_{-}\bigr )F +
  \left\{ I - \phi\bigl(T^\alpha \Lambda_{-}\bigr) \right\}F
  = F_{1} + F_{2},
$$
where $I$ denotes the identity operator. In view of
\eqref{LambdaMinusSymbolEquivalence},
\begin{align}
  \label{F1Spectrum}
  \bigabs{ \abs{\tau} - \abs{\xi} } &\le 4c^{2} T^{-\alpha}
  \quad \text{for} \quad (\tau,\xi) \in \supp \widehat{F_{1}}, \\
  \label{F2Spectrum}
  \bigabs{ \abs{\tau} - \abs{\xi} } &\ge T^{-\alpha}
  \quad \text{for} \quad (\tau,\xi) \in \supp \widehat{F_{2}}.
\end{align}
Now define
$$
  u = \chi(t) u_{0} + \chi(t/T) u_{1} + u_{2},
$$
where
\begin{align*}
  u_{0} &= \partial_{t} W(t)f + W(t) g \quad \text{with} \quad
  W(t) = \frac{\sin(t \sqrt{-\Delta})}{\sqrt{-\Delta}},
  \\
  u_{1} &= - \int_0^t W(t - t') F_{1}(t',\cdot) \, dt',
  \\
  \widehat u_{2} &= (\tau^2 - \abs{\xi}^2)^{-1} \widehat F_{2}.
\end{align*}
Observe that $F_1 \in H^{s-1,0} \subseteq L^1_{\mathrm{loc}}(\R,H^{s-1})$, so $u_1$ is
well-defined.
\begin{lemma} $u$ solves \eqref{CauchyProblem} on $[0,T] \times \R^n$.
\end{lemma}
\begin{proof} The only point which is not evident is that
$u_2 \init = \partial_t u_2 \init = 0$. This is clearly true
when $F \in \Schwartz$, since then $F_2 \in \Schwartz$, so $u_2$ is
necessarily given by Duhamel's formula. The general case then
follows by density, since clearly $F \mapsto u_2$
is linear and bounded from $H^{s-1,\theta-1}$ into $\scrH^{s,\theta}$,
and the latter space embeds in \eqref{GlobalEnergySpace}.
\end{proof}
Using hypothesis (iii) only, we shall prove the following lemmas.
\begin{lemma}\label{u0Lemma}
For all $(f,g) \in H^s \times H^{s-1}$,
$$
  \norm{\chi(t)u_{0}}_{\X^s} \le C \datanorm{f}{g}{s}.
$$
\end{lemma}
\begin{lemma}\label{u1Lemma}
For all $\varepsilon \ge 0$ and $F \in H^{s-1,\theta+\varepsilon-1}$,
$$
  \norm{\chi(t/T)u_{1}}_{\X^s} \le C T^\delta
  \spacetimenorm{F}{s-1}{\theta+\varepsilon-1},
$$
where
\begin{equation}\label{delta}
  \delta = 2 - \gamma - \alpha/2
  + \min \bigl\{ 0, \alpha(\theta + \varepsilon - 1) \bigr\},
\end{equation}
$\gamma$ is as in hypothesis (iii),
and $C$ is independent of $\varepsilon$ and $\alpha$.
\end{lemma}
Since by (i) we have \eqref{Yembedding}, and since
by \eqref{F2Spectrum} and (ii) we clearly have
\begin{equation}
  \label{u2Estimate}
  \norm{u_{2}}_{\X^s} \le C T^{\alpha \varepsilon} 
  \norm{F}_{\Y^{s,\varepsilon}},
\end{equation}
the theorem follows.

We turn to the proofs of Lemmas \ref{u0Lemma} and \ref{u1Lemma}.

Some notation: for $\gamma \in \R$, let $\D^\gamma$ be defined by
$$
  \widehat{\D^{\gamma}\chi}(\tau) = 
  (1 + \abs{\tau}^{2})^{\gamma/2} \widehat{\chi}(\tau).
$$
In what follows, $p \lesssim q$ means that $p \le C q$ for some
positive constant $C$ independent of $\alpha$ and $\varepsilon$.
\subsection{Proof of Lemma \ref{u1Lemma}}
We write $F_1 = F_{1,1} + F_{1,2}$, where
$\widehat{F_{1,1}}(\tau,\xi)$ and $\widehat{F_{1,2}}(\tau,\xi)$ are 
supported in the regions $\abs{\xi} \le T^{-\alpha}$ and $\abs{\xi} \ge 
T^{-\alpha}$ respectively. Let $u_{1,j}$ be defined as $u_{1}$,
but with $F_{1}$ replaced by $F_{1,j}$ for $j = 1, 2$.

The following lemma, proved in section \ref{ProofOfu1Characterization},
characterizes $u_{1,1}$ and $u_{1,2}$.
\begin{lemma}\label{u1Characterization}
Given $0 < T < 1$, let $u_{1,j}$, $j = 1, 2$ be defined as above.
There exist sequences $f_{j}^{+}, f_{j}^{-} \in H^{s}$ and
$g_{j} \in C([0,1],H^{s-1})$ such that
\begin{gather}
  \notag
  \begin{split}
    \supp \widehat{f_j^\pm} &\subseteq \{ \xi : \abs{\xi} \ge 
    T^{-\alpha} \}, \\
    \supp \widehat{g_j(\rho)} &\subseteq \{ \xi : \abs{\xi} \le 
    T^{-\alpha} \},
  \end{split}
  \\
  \label{u1CharacterizationA}
  \Sobnorm{f_{j}^{\pm}}{s},
  \sup_{0 \le \rho \le 1} \Sobnorm{g_{j}(\rho)}{s-1}
  \lesssim T^{\alpha(1/2-j)} \spacetimenorm{F_{1}}{s-1}{0}
\end{gather}
for $j = 1,2,\dots$, and
\begin{equation}\label{FirstFormulaForu12}
  u_{1,1} = \sum\nolimits_{1}, \qquad
  u_{1,2} = \sum\nolimits_{2} + E,
\end{equation}
where
\begin{align*}
  \sum\nolimits_{1} &= \sum_{j=1}^{\infty} \frac{t^{j+1}}{j!}
  \int_{0}^{1} e^{it(2\rho-1)\sqrt{-\Delta}} g_{j}(\rho) \, d\rho, \\
  \sum\nolimits_{2} &= \sum_{j=1}^{\infty} \frac{t^{j}}{j!} \left( e^{it\sqrt{-\Delta}} 
  f_{j}^{+} + e^{-it\sqrt{-\Delta}} f_{j}^{-} \right)
\end{align*}
and, given $\gamma \ge 1$,
\begin{multline}\label{EstimateForE}
  \sup_{\tau \in \R} \abs{ \Fourier \Lambda^{s} \Lambda_{-}^{\gamma}
  \{ \chi E \}(\tau,\xi) } \\
  \lesssim
  \left( \biginftynorm{\widehat{\D^{\gamma-1}\chi}}{} + 
  \biginftynorm{\widehat{\D^{\gamma}(t\chi)}}{} \right)
  \int \abs{ \Fourier \Lambda^{s-1} F_{1,2} (\lambda,\xi) } \, d\lambda
\end{multline}
for all $\xi \in \R^{n}$ and $\chi \in C_{c}^{\infty}(\R)$. Moreover,
\begin{equation}\label{SecondFormulaForu12}
  \widehat{\chi u_{1,2}}(\tau,\xi)
  = \frac{1}{2\pi}
  \int \int_{0}^{1} \! \int_{0}^{1}
  (1-\rho) \widehat{\chi}'' ( \tau - \mu )
  \widehat{F_{1,2}}(\lambda,\xi)
  \, d\sigma \, d\rho \, d\lambda,
\end{equation}
where
\begin{equation}\label{u1CharacterizationB}
  \mu = b + \sigma ( a - b ), \quad
  a = \abs{\xi} + \rho( \lambda - \abs{\xi} ), \quad
  b = - \abs{\xi} + \rho( \lambda + \abs{\xi} )
\end{equation}
and $\chi \in C_{c}^{\infty}(\R)$.
\end{lemma}
We need two more lemmas.
\begin{lemma}\label{CutOffLemma}
If $\sigma \in \R$, $\gamma \ge 0$, $f \in H^{\sigma}$ and $\chi \in 
C_{c}^{\infty}(\R)$, then
\begin{equation}\label{CutOffEstimate1}
  \Bdualnorm{ \Fourier \Lambda^{\sigma} \Lambda_{-}^{\gamma}
  \{ \chi(t) e^{\pm i t \sqrt{-\Delta}} f \}(\tau,\xi) }
  \lesssim \biginftynorm{\widehat{\D^{\gamma} \chi}}{} \Sobnorm{f}{\sigma}.
\end{equation}
If $-1 \le \rho \le 1$, $g \in H^{\sigma}$ and $\supp \widehat{g} \subseteq 
\{ \xi : \abs{\xi} \le r \}$, then
\begin{equation}\label{CutOffEstimate2}
  \Bdualnorm{ \Fourier \Lambda^{\sigma} \Lambda_{-}^{\gamma}
  \{ \chi(t) e^{i \rho t \sqrt{-\Delta}} g \}(\tau,\xi) }
  \lesssim \left( \biginftynorm{\widehat{\D^{\gamma} \chi}}{}
  + r^{\gamma} \inftynorm{\widehat{\chi}}{} \right)
  \Sobnorm{g}{\sigma}.
\end{equation}
\end{lemma}
\begin{proof}
This is a triviality. To prove \eqref{CutOffEstimate1}, simply note 
that the Fourier transform of $\chi(t) e^{\pm it\sqrt{-\Delta}} f$ is
$\widehat{\chi}(\tau \mp \abs{\xi}) \widehat{f}(\xi)$.
To prove, \eqref{CutOffEstimate2} observe that
the Fourier transform of $\chi(t) e^{i\rho t \sqrt{-\Delta}}g$ is
$\widehat{\chi}( \tau - \rho \abs{\xi} ) \widehat{g}(\xi)$
and that
$$
  \hypwt{\tau}{\xi} \le \bigl\vert \tau - \rho \abs{\xi} 
  \bigr\vert + ( 1 - \abs{\rho} ) \abs{\xi}
  \le \bigl\vert \tau - \rho \abs{\xi} \bigr\vert + r
$$
for $\xi \in \supp \widehat{g}$ and $-1 \le \rho \le 1$.
\end{proof}
\begin{lemma}\label{CutOffLemmaCorollary}
Assume that the inequality
$$
  \norm{u}_{\X^{s}} \lesssim \bignorm{ \Fourier \Lambda^{s-1} \Lambda_{+}
  \Lambda_{-}^{\gamma} u (\tau,\xi) }_{L_{\xi}^{2}(L_{\tau}^{\infty})}
$$
of hypothesis (iii) of Theorem \ref{BasicConditionsTheorem} holds
for some $\gamma \in \R$. Then for all $\chi \in C_{c}^{\infty}(\R)$ 
and $f \in H^{s}$,
\begin{equation}\label{CutOffLemmaEstimate1}
  \norm{ \chi(t) e^{\pm it\sqrt{-\Delta}} f}_{\X^s}
  \lesssim \biginftynorm{\widehat{\D^{\gamma}\chi}}{} \Sobnorm{f}{s}
  + \biginftynorm{\widehat{\D^{\gamma+1}\chi}}{} \Sobnorm{f}{s-1}.
\end{equation}
Moreover, if $-1 \le \rho \le 1$, $g \in H^{s-1}$ and
$\supp \widehat{g} \subseteq \{ \xi : \abs{\xi} \le r \}$, where
$r \ge 1$, then
\begin{equation}\label{CutOffLemmaEstimate4}
  \norm{ \chi(t) e^{i\rho t \sqrt{-\Delta}} g }_{\X^{s}}
  \lesssim \left( r^{\gamma+1}\inftynorm{\widehat{\chi}}{}
  + r \biginftynorm{\widehat{\D^{\gamma}\chi}}{}
  + \biginftynorm{\widehat{\D^{\gamma+1}\chi}}{}\right) 
  \Sobnorm{g}{s-1}.
\end{equation}
\end{lemma}
\begin{proof}
Observe that
\begin{multline*}
  \Bdualnorm{\Fourier \Lambda^{s-1} \Lambda_{+} \Lambda_{-}^{\gamma} 
  u(\tau,\xi)} \\
  \lesssim
  \Bdualnorm{\Fourier \Lambda^{s} \Lambda_{-}^{\gamma} 
  u(\tau,\xi)}
  + \bigBdualnorm{\Fourier \Lambda^{s-1} \Lambda_{-}^{\gamma+1} 
  u(\tau,\xi)}
\end{multline*}
for any $u(t,x)$. Thus, \eqref{CutOffLemmaEstimate1} and
\eqref{CutOffLemmaEstimate4} follow immediately from
\eqref{CutOffEstimate1} and \eqref{CutOffEstimate2} of Lemma 
\ref{CutOffLemma}, respectively.
\end{proof}
We are now in a position to finish the proof of Lemma \ref{u1Lemma}.
\paragraph{Estimate for $u_{1,1}$.}
By \eqref{F1Spectrum},
\begin{equation}\label{F1Estimate}
  \spacetimenorm{F_{1}}{s-1}{0} \lesssim T^{\min\{ 0, \alpha(\theta + 
  \varepsilon - 1)\}} \spacetimenorm{F}{s-1}{\theta+\varepsilon-1}.
\end{equation}
It is easily checked that
\begin{equation}\label{ScalingEstimate}
  \inftynorm{\Fourier \D^{r} \bigl\{ t^{j}\chi(t/T) \bigr\}}{}
  \le T^{1-r+j}
  \biginftynorm{\Fourier \D^{r}(t^{j}\chi)}{}
\end{equation}
for $r \ge 0$, $j \ge 0 $ and $0 < T < 1$. Combining
\eqref{CutOffLemmaEstimate4} of Lemma \ref{CutOffLemmaCorollary}
with \eqref{u1CharacterizationA}, \eqref{F1Estimate}
and \eqref{ScalingEstimate} yields
$$
  \norm{\chi(t/T)\sum\nolimits_{1}}_{\X^{s}}
  \le C_{\chi,\gamma} T^\delta
  \spacetimenorm{F_{1}}{s-1}{\theta+\varepsilon-1},
$$
where $\delta$ is given by \eqref{delta} and
$$
  C_{\chi,\gamma} \lesssim \sum_{j=1}^{\infty}
  \frac{1}{j!}
  \biginftynorm{ \Fourier{ \D^{\gamma+1}( t^{j+1} \chi ) }}{}
  \lesssim \sum_{j=1}^{\infty}
  \frac{1}{j!} \left( \Lpnorm{ t^{j+1} \chi}{1}{}
  + \Lpnorm{ (t^{j+1} \chi)'''}{1}{} \right).
$$
Since $\chi$ is compactly supported, $C_{\chi,\gamma} < \infty$.
\paragraph{Estimate for $u_{1,2}$.}
In view of hypothesis (iii) of the theorem,
$$
  \norm{ \chi(t/T)u_{1,2} }_{\X^{s}}
  \lesssim A + B,
$$
where
\begin{align*}
  A &= \norm{ \Fourier \Lambda^{s} \Lambda_{-}^{\gamma}
  \{ \chi(t/T) u_{1,2}\}(\tau,\xi) }_{L_{\xi}^{2}(L_{\tau}^{\infty})
  (\abs{\tau} \le C \abs{\xi})}, \\
  B &= \norm{ \Fourier \Lambda^{s-1} \Lambda_{-}^{\gamma+1}
  \{ \chi(t/T) u_{1,2}\}(\tau,\xi) }_{L_{\xi}^{2}(L_{\tau}^{\infty})
  (\abs{\tau} \ge C \abs{\xi})}
\end{align*}
and $C > 1$ is a constant which will be specified later.

To estimate $A$ we use the expression \eqref{FirstFormulaForu12} for
$u_{1,2}$. By \eqref{CutOffEstimate1} of Lemma \ref{CutOffLemma},
as well as \eqref{u1CharacterizationA}, \eqref{F1Estimate}
and \eqref{ScalingEstimate},
$$
  \Bdualnorm{ \Fourier \Lambda^{s} \Lambda_{-}^{\gamma}
  \left\{ \chi(t/T) \sum\nolimits_{2} \right\}(\tau,\xi) }
  \lesssim C_{\chi,\gamma} T^\delta
  \spacetimenorm{F_{1}}{s-1}{\theta+\varepsilon-1},
$$
where
$$
  C_{\chi,\gamma} \lesssim \sum_{j=1}^{\infty}
  \frac{1}{j!} \inftynorm{\Fourier\D^{\gamma}(t^{j}\chi)}{} < \infty.
$$
Next, by \eqref{EstimateForE} and \eqref{F1Spectrum},
\begin{multline*}
  \Bdualnorm{ \Fourier \Lambda^{s} \Lambda_{-}^{\gamma}
  \{ \chi(t/T) E \}(\tau,\xi) }
  \\
  \lesssim
  \left( \biginftynorm{\widehat{\D^{\gamma-1}\chi}}{} + 
  \biginftynorm{\widehat{\D^{\gamma}(t\chi)}}{} \right)
  T^{-\alpha/2} \spacetimenorm{F_{1}}{s-1}{0},
\end{multline*}
which combined with \eqref{F1Estimate} and \eqref{ScalingEstimate} 
yields
$$
  \Bdualnorm{ \Fourier \Lambda^{s} \Lambda_{-}^{\gamma}
  \{ \chi(t/T) E \}(\tau,\xi) }
  \lesssim C_{\chi,\gamma} T^\delta
  \spacetimenorm{F_{1}}{s-1}{\theta+\varepsilon-1},
$$
where
$$
  C_{\chi,\gamma} \lesssim
  \biginftynorm{\widehat{\D^{\gamma-1}\chi}}{} + 
  \biginftynorm{\widehat{\D^{\gamma}(t\chi)}}{}.
$$

For $B$ we use the expression \eqref{SecondFormulaForu12}
for $\widehat{ \chi u_{1,2} }(\tau,\xi)$. Recall that
$\abs{\tau} \ge C \abs{\xi}$ in the definition of $B$, where
$C$ is yet to be determined.
We claim that if we take $C = 8(1+c^{2})$, where $c$ is the constant
appearing in \eqref{LambdaMinusSymbolEquivalence} and \eqref{F1Spectrum},
then with $\mu$ given by \eqref{u1CharacterizationB},
$$
  \abs{\tau} \ge C \abs{\xi} \implies
  \hypwt{\tau}{\xi} \sim \abs{\tau} \sim  \abs{ \tau - \mu }
$$
for all $\lambda$ such that $(\lambda,\xi) \in \supp 
\widehat{F_{1,2}}$, and all $0 \le \rho, \sigma \le 1$. Indeed,
$$
  \abs{\mu} \le \abs{b} + \abs{a-b}
  \le 4 \abs{\xi} + \hypwt{\lambda}{\xi},
$$
and we have
$$
  \hypwt{\lambda}{\xi} \le 4c^{2} T^{-\alpha}
  \le 4c^{2} \abs{\xi}
$$
for $(\lambda,\xi) \in \supp \widehat{F_{1,2}}$. Thus,
$$
  B \lesssim \biginftynorm{ \Fourier \D^{\gamma + 1} \{ 
  t^{2}\chi(t/T) \} }{}
  \int \abs{ \Fourier \Lambda^{s-1} F_{1,2} (\lambda,\xi) } \, d\lambda,
$$
and using \eqref{F1Estimate} and \eqref{ScalingEstimate}, we get
$$
  B
  \lesssim C_{\chi,\gamma} T^\delta
  \spacetimenorm{F_{1}}{s-1}{\theta+\varepsilon-1},
$$
where $C_{\chi,\gamma} \lesssim \inftynorm{ \Fourier \D^{\gamma + 1} ( 
t^{2}\chi ) }{}$.
\subsection{Proof of Lemma \ref{u0Lemma}}
We will show, using hypothesis (iii) of Theorem \ref{BasicConditionsTheorem},
that for all $\chi \in C_{c}^{\infty}(\R)$ 
and $(f,g) \in H^{s} \times H^{s-1}$,
\begin{equation}\label{CutOffLemmaEstimate2}
  \norm{ \chi(t) \partial_t W(t) f}_{\X^s}
  \lesssim \biginftynorm{\widehat{\D^{\gamma}\chi}}{} \Sobnorm{f}{s}
  + \biginftynorm{\widehat{\D^{\gamma+1}\chi}}{} \Sobnorm{f}{s-1}
\end{equation}
and
\begin{multline}\label{CutOffLemmaEstimate3}
  \norm{ \chi(t) W(t) g}_{\X^s}
  \lesssim \left( \biginftynorm{\widehat{\D^{\gamma}\chi}}{}
  + \biginftynorm{\widehat{\D^{\gamma+1}(t\chi)}}{} \right)
  \Sobnorm{g}{s-1}
  \\
  + \biginftynorm{\widehat{\D^{\gamma+1}\chi}}{} \Sobnorm{g}{s-2}.
\end{multline}
Clearly this implies Lemma \ref{u0Lemma}.

We apply Lemma \ref{CutOffLemmaCorollary}. Evidently, \eqref{CutOffLemmaEstimate1}
implies \eqref{CutOffLemmaEstimate2}, so it remains to prove
\eqref{CutOffLemmaEstimate3}. For this, we split $g = g_{1} + g_{2}$,
where $\widehat{g_{1}}$ is supported in the region $\abs{\xi} < 1$
and $\widehat{g_{2}}$ is supported in $\abs{\xi} \ge 1$.
Since
$$
  (-\Delta)^{-\half} \sin(t\sqrt{-\Delta}) = t \int_{0}^{1} e^{it(2\rho-1)\sqrt{-\Delta}} \, d\rho,
$$
we have
$$
  \chi(t) (-\Delta)^{-\half} \sin(t\sqrt{-\Delta}) g_{1} = \int_{0}^{1} t\chi(t)
  e^{it(2\rho-1)\sqrt{-\Delta}} g_{1} \, d\rho.
$$
By \eqref{CutOffLemmaEstimate4},
$$
  \norm{ t\chi(t) e^{it(2\rho-1)\sqrt{-\Delta}} g_{1} }_{\X^{s}}
  \lesssim \biginftynorm{\Fourier{\D^{\gamma+1}(t\chi)}}{}
  \Sobnorm{g_{1}}{s-1}
$$
for $0 \le \rho \le 1$, and it follows that
$$
  \norm{ \chi(t) W(t) g_{1}}_{\X^{s}}
  \le \biginftynorm{\widehat{\D^{\gamma+1}(t\chi)}}{}
  \Sobnorm{g_{1}}{s-1}.
$$
This proves \eqref{CutOffLemmaEstimate3} with $g$ replaced by
its low frequency part $g_{1}$.

Since $\bigSobnorm{(-\Delta)^{-\half} g_{2}}{s} \le 2 \Sobnorm{g}{s-1}$, the estimate
\eqref{CutOffLemmaEstimate3} with $g$ replaced by $g_{2}$ follows 
immediately from \eqref{CutOffLemmaEstimate1}.
\subsection{Proof of Lemma \ref{u1Characterization}}\label{ProofOfu1Characterization}
We write
$u_{1} = - (2i)^{-1} (-\Delta)^{-\half}(v_{+} - v_{-})$, where
$$
  v_{\pm}(t) = \int_{0}^{t} e^{\pm i(t-t')\sqrt{-\Delta}} F_{1}(t',\cdot) \, dt'.
$$
As in \cite[Lemma 4.4]{Kl-Ma1},
\begin{align}
  \label{FirstFormulaForv}
  \widehat{v_{\pm}(t)}(\xi)
  &= \frac{e^{\pm it\abs{\xi}}}{2\pi} \int
  \frac{e^{it(\tau \mp \abs{\xi})} - 1}{i(\tau \mp \abs{\xi})}
  \widehat{F_{1}}(\tau,\xi) \, d\tau.
  \\
  \label{SecondFormulaForv}
  & = \frac{e^{\pm it\abs{\xi}}}{2\pi}
  \sum_{j=1}^{\infty} \frac{t^{j}}{j!}
  \int i^{j-1} (\tau \mp \abs{\xi})^{j-1}
  \widehat{F_{1}}(\tau,\xi) \, d\tau.
\end{align}
\paragraph{Formula for $u_{1,1}$.}
By \eqref{SecondFormulaForv}, we have
\begin{align*}
  \widehat{u_{1,1}(t)}(\xi)
  &= \frac{1}{4\pi} \sum_{j=1}^{\infty} \frac{t^{j}}{j!}
  \int i^{j} \abs{\xi}^{-1} \bigl( \beta_{j,\tau}(\abs{\xi})
  - \beta_{j,\tau}(-\abs{\xi}) \bigr) \widehat{F_{1,1}}(\tau,\xi) \, d\tau
  \\
  &= \frac{1}{2\pi} \sum_{j=1}^{\infty} \frac{t^{j}}{j!}
  \int_{0}^{1} \! \int i^{j} \beta_{j,\tau}' 
  \bigl( (2\rho - 1)\abs{\xi} \bigr)
  \widehat{F_{1,1}}(\tau,\xi) \, d\tau \, d\rho,
\end{align*}
where $\beta_{j,\tau}(r) = e^{itr}( \tau - r )^{j-1}$. Since
$$
  \beta_{j,\tau}'(r) = it e^{itr} ( \tau - r )^{j-1}
  - e^{itr} (j-1) ( \tau - r )^{j-2},
$$
where the second term only occurs for $j \ge 2$, we get
$$
  u_{1,1}(t)
  = \frac{1}{2\pi} \sum_{j=1}^{\infty} \frac{t^{j+1}}{j!}
  \int_{0}^{1} i^{j+1} \Bigl( 1 - \frac{j}{j+1} \Bigr)
  e^{it(2\rho-1)\sqrt{-\Delta}} k_{j}(\rho) \, d\rho,
$$
where $k_{j}(\rho)$ is given by
$$
  \widehat{k_{j}(\rho)}(\xi)
  = \int \bigl( \tau - (2\rho-1)\abs{\xi} \bigr)^{j-1} 
  \widehat{F_{1,1}}(\tau,\xi) \, d\tau.
$$
If we set
$$
  g_{j} = \frac{i^{j+1}}{2\pi} \left( 1 - \frac{j}{j+1} \right) k_{j},
$$
then
$$
   u_{1,1}(t) = \sum_{j=1}^{\infty} \frac{t^{j+1}}{j!}
  \int_{0}^{1} e^{it(2\rho-1)\sqrt{-\Delta}} g_{j}(\rho) \, d\rho.
$$
Notice that
$$
  \bigl\vert \tau - (2\rho-1) \abs{\xi} \bigr\vert
  \le \abs{\tau} + \abs{\xi} \le \hypwt{\tau}{\xi} + 2\abs{\xi}
$$
for $0 \le \rho \le 1$. But if $(\tau,\xi) \in \supp \widehat{F_{1,1}}$,
then $\abs{\xi} \le T^{-\alpha}$, and by \eqref{F1Spectrum} we also have
$\hypwt{\tau}{\xi} \lesssim T^{-\alpha}$. It follows that
$$
  \bigabs{\widehat{g_{j}(\rho)}(\xi)}
  \lesssim T^{\alpha(1-j)} \int \bigabs{\widehat{F_{1,1}}(\tau,\xi)} \, d\tau
  \lesssim T^{\alpha(1/2-j)}
  \left( \int \bigabs{\widehat{F_{1,1}}(\tau,\xi)}^{2} \, d\tau 
  \right)^{\half},
$$
whence $\Sobnorm{g_{j}(\rho)}{s-1} \lesssim T^{\alpha(1/2-j)} 
\spacetimenorm{F_{1}}{s-1}{0}$ for $0 \le \rho \le 1$.
\paragraph{First formula for $u_{1,2}$.}
By \eqref{FirstFormulaForv} and \eqref{SecondFormulaForv},
$$
  u_{1,2}(t) = \sum_{j=1}^{\infty} \frac{t^{j}}{j!} \left( e^{it\sqrt{-\Delta}} 
  f_{j}^{+} + e^{-it\sqrt{-\Delta}} f_{j}^{-} \right) + E_{+}(t) + E_{-}(t),
$$
where
\begin{align*}
  \widehat{f_{j}^{+}}(\xi) &= (4 \pi \abs{\xi})^{-1}
  \int_{0}^{\infty} i^{j} (\abs{\tau} - \abs{\xi})^{j-1}
  \widehat{F_{1,2}}(\tau,\xi) \, d\tau, \\
  \widehat{f_{j}^{-}}(\xi) &= - (4 \pi \abs{\xi})^{-1}
  \int_{-\infty}^{0} i^{j} (\abs{\xi} - \abs{\tau})^{j-1}
  \widehat{F_{1,2}}(\tau,\xi) \, d\tau,\\
  \widehat{E_{+}(t)}(\xi) &= - \frac{1}{4\pi\abs{\xi}} \int_{-\infty}^{0}
  \frac{e^{it\tau} - e^{it\abs{\xi}}}{\abs{\tau} + 
  \abs{\xi}} \widehat{F_{1,2}}(\tau,\xi) \, d\tau, \\
  \widehat{E_{-}(t)}(\xi) &= - \frac{1}{4\pi\abs{\xi}} \int_{0}^{\infty}
  \frac{e^{it\tau} - e^{-it\abs{\xi}}}{\abs{\tau} + 
  \abs{\xi}} \widehat{F_{1,2}}(\tau,\xi) \, d\tau.
\end{align*}
It follows easily from \eqref{F1Spectrum} that
$\Sobnorm{f_{j}^{\pm}}{s} \lesssim T^{\alpha(1/2-j)} \spacetimenorm{F_{1}}{s-1}{0}$.

Next, observe that
$$
  \widehat{\chi E_{+}}(\tau,\xi)
  = - \frac{1}{4\pi\abs{\xi}} \int_{-\infty}^{0}
  \frac{ \widehat{\chi}(\tau-\lambda) - \widehat{\chi}(\tau-\abs{\xi})}
  {\abs{\lambda} + \abs{\xi}} \widehat{F_{1,2}}(\lambda,\xi) \, 
  d\lambda.
$$
Since
$$
  \frac{ \widehat{\chi}(\tau-\lambda) -
  \widehat{\chi}(\tau-\abs{\xi})}{\abs{\lambda} + \abs{\xi}}
  = \int_{0}^{1} \widehat{\chi}'
  \bigl( \tau - \abs{\xi} + \rho ( \abs{\lambda} + \abs{\xi} )
  \bigr) \, d\rho
$$
for $\lambda < 0$, it is easy to see, by considering separately
the two cases $\hypwt{\tau}{\xi} \le 2(\abs{\lambda} + \abs{\xi})$
and $\hypwt{\tau}{\xi} > 2(\abs{\lambda} + \abs{\xi})$, that for a
given $\gamma \ge 1$,
$$
  \bigl( 1 + \hypwt{\tau}{\xi} \bigr)^{\gamma}
  \frac{ \bigabs{\widehat{\chi}(\tau-\lambda)
  - \widehat{\chi}(\tau-\abs{\xi}) } }{ \abs{\lambda} + \abs{\xi} }
  \lesssim \biginftynorm{\widehat{\D^{\gamma-1}\chi}}{} + 
  \biginftynorm{\widehat{\D^{\gamma}(t\chi)}}{}
$$
for all $\tau \in \R$, $\xi \in \R^{n}$ and $\lambda < 0$.
We conclude that
\begin{multline*}
  \abs{ \Fourier \Lambda^{s} \Lambda_{-}^{\gamma}
  \{ \chi E_{+} \}(\tau,\xi) } \\
  \lesssim
  \left( \biginftynorm{\widehat{\D^{\gamma-1}\chi}}{} + 
  \biginftynorm{\widehat{\D^{\gamma}(t\chi)}}{} \right)
  \int \abs{ \Fourier \Lambda^{s-1} F_{1,2} (\lambda,\xi) } \, d\lambda.
\end{multline*}
The same estimate holds for $E_{-}$, so
$E = E_{+} + E_{-}$ satisfies \eqref{EstimateForE}.
\paragraph{Second formula for $u_{1,2}$.}
Using \eqref{FirstFormulaForv}, we write
$$
  \widehat{u_{1,2}(t)}(\xi) = \frac{1}{4\pi\abs{\xi}}
  \int \left\{ \frac{ e^{it\tau} - e^{it\abs{\xi}}}{\tau - \abs{\xi}}
  - \frac{ e^{it\tau} - e^{-it\abs{\xi}}}{\tau + \abs{\xi}}
  \right\} \widehat{F_{1,2}}(\tau,\xi) \, d\tau.
$$
Thus,
\begin{align*}
  &\widehat{\chi u_{1,2}}(\tau,\xi) \\
  &= \frac{1}{4\pi\abs{\xi}}
  \int \left\{ \frac{ \widehat{\chi}(\tau - \lambda)
  - \widehat{\chi}(\tau - \abs{\xi})}{\lambda - \abs{\xi}}
  - \frac{ \widehat{\chi}(\tau - \lambda)
  - \widehat{\chi}(\tau + \abs{\xi})}{\lambda + \abs{\xi}}
  \right\} \widehat{F_{1,2}}(\lambda,\xi) \, d\lambda \\
  &= - \frac{1}{4\pi\abs{\xi}}
  \int \! \int_{0}^{1} \bigl\{
  \widehat{\chi}' (\tau - a)
  - \widehat{\chi}' (\tau - b) \bigr\}
  \widehat{F_{1,2}}(\lambda,\xi) \, d\rho \, d\lambda \\
  &= \frac{1}{2\pi}
  \int \int_{0}^{1} \! \int_{0}^{1}
  (1-\rho) \widehat{\chi}'' \bigl( \tau - b + \sigma(b-a) \bigr)
  \widehat{F_{1,2}}(\lambda,\xi)
  \, d\sigma \, d\rho \, d\lambda,
\end{align*}
where $a = \abs{\xi} + \rho( \lambda - \abs{\xi})$ and
$b = - \abs{\xi} + \rho( \lambda + \abs{\xi})$.
\section{Proof of Theorem \ref{AbstractWPTheorem}}\label{ProofOfAbstractWPTheorem}
We may assume that the function $A$ in 
\eqref{NonlinearEstimate} and \eqref{NonlinearDifferenceEstimate} is 
increasing, and that $C_{T}$ is an increasing and continuous function 
of $T$.
\subsection{Local existence}\label{Existence}
Given $(f,g)$ and $0 < T < 1$,
let $u_{0}$ be the solution of the homogeneous wave equation with 
initial data $(f,g)$, and set
$$
  u_{j} = u_{0} + W \mathcal N(u_{j-1}), \quad j = 1,2,\dots
$$
By (S2), there is a constant $C$ such that
\begin{equation}\label{HomogeneousEstimate}
  \norm{u_{0}}_{\X^{s}_{T}} \le C \bigl( \Sobnorm{f}{s} + 
  \Sobnorm{g}{s-1} \bigr).
\end{equation}
Combined with \eqref{NonlinearEstimate} this gives
$$
  \norm{u_{j}}_{\X^{s}_{T}} \le R/2 + C_{T} A \bigl( 
  \norm{u_{j-1}}_{\X^{s}_{T}} \bigr) \norm{u_{j-1}}_{\X^{s}_{T}}
$$
for $j \ge 1$, where $R$ is twice the right hand side of 
\eqref{HomogeneousEstimate}. By \eqref{CTDecay},
we may choose $T$ so small that $2 C_{T} A(R) \le 1$.
Since $A$ is increasing, it now follows by
induction that $\norm{u_{j}}_{\X^{s}_{T}} \le R$ for all $j$.
It then follows by \eqref{NonlinearDifferenceEstimate} that
\begin{equation}\label{CauchyCondition}
  \norm{u_{j+1}-u_{j}}_{\X^{s}_{T}}
  \le \half \norm{u_{j}-u_{j-1}}_{\X^{s}_{T}}
\end{equation}
for $j \ge 1$. Thus, $(u_{j})$ is a Cauchy sequence in $\X^{s}_{T}$,
and we let $u$ be its limit. In view of 
\eqref{NonlinearDifferenceEstimate}, $\mathcal N(u_{j}) \to \mathcal 
N(u)$ in the sense of distributions on $S_T = (0,T) \times \R^{n}$. Since
$\square u_{j} = \mathcal N(u_{j-1})$ on $S_{T}$
with initial data $(f,g)$, by passing to the limit we conclude that
$\square u = \mathcal N(u)$ on $S_{T}$ with the same data.
\subsection{Uniqueness}\label{Uniqueness}
Assume that $T > 0$ and $u,u' \in \X^s_T$ are two solutions of
\eqref{GenericSystem} on $S_T$ with the same
initial data $(f,g)$. It suffices to prove that the set
$$
  E = \{ t \in [0,T] : \text{$u(\rho) = u'(\rho)$ for all $\rho \in [0,t]$} \}
$$
is open in $[0,T]$, since $E$ is obviously closed and nonempty.

Assume that $t \in E$, $t < T$. By (S4) we 
may consider $u$ and $u'$ to be elements of $\X^{s}_{[t,T]}$, and by 
(N3) they are both solutions of \eqref{GenericSystem} on $(t,T) \times \R^{n}$ 
with the same initial data at time $t$ (since $t \in E$). Next, by 
(S3) and (N2), $\tau_t u, \tau_t u' \in \X^{s}_{T-t}$ solve
\eqref{GenericSystem} on $(0,T-t) \times \R^{n}$ with identical
initial data at time $0$.

By the above, it suffices to prove that $\varepsilon \in E$ for some 
arbitrarily small $\varepsilon > 0$. But by \eqref{NonlinearDifferenceEstimate},
$$
  \norm{u-u'}_{\X^s_{\varepsilon}} \le B(\varepsilon) 
  \norm{u-u'}_{\X^s_{\varepsilon}},
$$
where $B(\varepsilon) = C_{\varepsilon}
A\bigl( \max\{ \norm{u}_{X^s_\varepsilon}, \norm{u'}_{X^s_\varepsilon} 
\} \bigr)$, and in view of \eqref{CTDecay} and property (S4),
$\lim_{\varepsilon \to 0^{+}} B(\varepsilon) = 0$.
\subsection{Continuous dependence on initial data}\label{ContDep}
\paragraph{Step 1.}
We prove that (III) follows from a weaker condition.
We denote by $u(f,g) \in \X^{s}_{T}$ the solution obtained in section
\ref{Existence}. Recall that $T = T(f,g) > 0$ is continuous, and
\begin{equation}\label{SolutionBound}
  C_{T(f,g)} A \bigl( \norm{u(f,g)}_{\X^{s}_{T(f,g)}} \bigr) \le \half.
\end{equation}
We claim that (III) follows from
\textit{
\begin{itemize}
  \item[(III$\,'$)]
  The map $(f,g) \mapsto u(f,g)$ is Lipschitz, in the sense that
  \begin{equation}\label{Lipschitz}
    \norm{ u(f,g) - u(f',g') }_{\X^s_{T}}
    \lesssim \Sobnorm{f - f'}{s} + \Sobnorm{g - g'}{s-1}
  \end{equation}
  for all initial data pairs $(f,g)$ and $(f',g')$ in $H^{s} \times 
  H^{s-1}$, where $T = \min\{ T(f,g), T(f',g') \}$.
\end{itemize}
}
With hypotheses as in (III), set
$$
  T_{*} = \inf_{0 \le t \le T} T\bigl( u(t),\partial_{t} u(t) \bigr).
$$
In view of (S1) and the continuity of $T$, $T_{*} > 0$.
Pick $0 < \varepsilon < T_{*}/2$ such that
$T = M \varepsilon$ for some integer $M$, and set $t_{j} = j 
\varepsilon$, $f_{j} = u(t_{j})$ and $g_{j} = \partial_{t} u(t_{j})$.

Assume that (III$'$) holds. In view of (S3) and (N2), for $j = 
0,1,\dots,M-2$ there is a ball $B_{j}$ in $H^{s} \times H^{s-1}$,
centered at $(f_{j}, g_{j})$, with the following property:
for all $(\phi_{j}, \psi_{j}) \in B_{j}$ there exists
$v_{j} \in \X^{s}_{[t_{j},t_{j+2}]}$ which solves
\eqref{GenericSystem} on $(t_{j},t_{j+2}) \times \R^{n}$
with initial data $(\phi_{j},\psi_{j})$ at time $t_{j}$, and satisfies
\begin{equation}\label{ContinuityA}
  \norm{u - v_{j}}_{\X^{s}_{[t_{j},t_{j+2}]}}
  \lesssim \Sobnorm{f_{j} - \phi_{j}}{s}
  + \Sobnorm{g_{j} - \psi_{j}}{s-1}.
\end{equation}
By (S1), this implies
\begin{multline}\label{ContinuityB}
  \Sobnorm{f_{j+1} - \phi_{j+1}}{s}
  + \Sobnorm{g_{j+1} - \psi_{j+1}}{s-1}
  \\
  \le C \bigl( \Sobnorm{f_{j} - \phi_{j}}{s}
  + \Sobnorm{g_{j} - \psi_{j}}{s-1} \bigl),
\end{multline}
where we have set $\phi_{j+1} = v(t_{j+1})$ and $\psi_{j+1} = \partial_{t} v(t_{j+1})$.

Thus, if we make $B_{M-3}$ so small that $C B_{M-3} \subseteq
B_{M-2}$, and then make $B_{M-4}$ so small that $C B_{M-4} 
\subseteq B_{M-3}$ etc., we find that if we start with data
$(\phi_{0},\psi_{0}) \in B_{0}$, then $v_{0}$ exists and
$$
  (\phi_{1},\psi_{1}) =
  \bigl(v_{0}(t_{1}), \partial_{t} u_{0}(t_{1}) \bigr) \in B_{1},
$$
so $v_{1}$ exists, and so on.

By translation invariance and uniqueness, the different
$v_{j}$ agree on the intersection of their domains, so by
(S5) we get a solution $v \in \X^{s}_{T}$ of \eqref{GenericSystem} on
$(0,T) \times \R^{n}$ with data $(\phi_{0},\psi_{0})$, and
$\norm{u - v}_{\X^{s}_{T}} \lesssim \sum_{j=0}^{M-2}
\norm{u - v_{j}}_{\X^{s}_{[t_{j},t_{j+2}]}}$.
But in view of \eqref{ContinuityA} and \eqref{ContinuityB},
$$
  \norm{u - v_{j}}_{\X^{s}_{[t_{j},t_{j+2}]}}
  \lesssim \Sobnorm{f-\phi_{0}}{s} + \Sobnorm{g-\psi_{0}}{s-1}.
$$
\paragraph{Step 2.}
We prove (III$'$). In view of \eqref{HomogeneousEstimate} and 
\eqref{NonlinearDifferenceEstimate}, it suffices to prove
\begin{equation}\label{ContinuousDependenceEstimate}
  C_{T} A \bigl( \max\{ \norm{u(f,g)}_{\X^s_{T}},
  \norm{u(f',g')}_{\X^s_{T}} \} \bigr)
  \le \half,
\end{equation}
where $T = \min\{ T(f,g), T(f',g') \}$.
But by \eqref{SolutionBound} and (S4),
$$
  A \bigl( \max\{ \norm{u(f,g)}_{\X^s_{T}},
  \norm{u(f',g')}_{\X^s_{T}} \} \bigr)
  \le \half \left( \min\{ C_{T(f,g)}, C_{T(f',g')} \} \right)^{-1}.
$$
Since we assume that $C_{T}$ is increasing in $T$,
\eqref{ContinuousDependenceEstimate} follows.
\subsection{Persistence of higher regularity}\label{Persistence}
We prove the assertion in remark (2) following Theorem
\ref{AbstractWPTheorem}.
\paragraph{Step 1.}
We show that (V) follows from
\textit{
\begin{itemize}
  \item[(I$\,'$)]
  Let $\sigma \ge s$. For all $(f,g) \in H^{\sigma} \times 
  H^{\sigma-1}$ there exist a $T > 0$ and a $u \in \X^s_{T} \cap 
  C\bigl( [0,T], H^{\sigma} \bigr) \cap
  C^{1} \bigl( [0,T], H^{\sigma-1} \bigr)$ which solves \eqref{GenericSystem} on
  $S_{T} = (0,T) \times \R^{n}$ with initial data $(f,g)$.
  Moreover, $T$ can be chosen to depend continuously on 
  $\Sobnorm{f}{s} + \Sobnorm{g}{s-1}$.
\end{itemize}
}
Assume that (I$'$) holds for a fixed $\sigma \ge s$, and let us denote
the existence time by $T(f,g)$. This function may depend on $s$ and $\sigma$,
but these are fixed quantities.

With hypotheses as in (V), set
$$
  T_{*} = \inf_{0 \le t \le T} T\bigl( u(t),\partial_{t} u(t) \bigr),
$$
and choose $0 < \varepsilon < T_{*}/2$ so that $T = M \varepsilon$ for
some integer $M$. Then set $t_{j} = j \varepsilon$, $f_{j} = u(t_{j})$
and $g_{j} = \partial_{t} u(t_{j})$.
By (I$'$) and translation invariance, for $j = 0,1,\dots,M-2$
there exists
$$
  u_{j} \in \X^{s}_{[t_{j},t_{j+2}]} \cap 
  C\bigl( [t_{j},t_{j+2}], H^{\sigma} \bigr) \cap
  C^{1} \bigl( [t_{j},t_{j+2}], H^{\sigma-1} \bigr)
$$
which solves \eqref{GenericSystem} 
on $(t_{j},t_{j+2}) \times \R^{n}$ with initial data $(f_{j},g_{j})$.
By uniqueness, each $u_{j}$
agrees with $u$ on $[t_{j},t_{j+2}]$. We conclude that 
$$
   u \in C\bigl( [0,T],H^{\sigma} \bigr) \cap C^{1}\bigl( 
   [0,T],H^{\sigma-1} \bigr).
$$
\paragraph{Step 2.}
We prove (I$'$). If we fix $\sigma \ge s$, we may assume that
the function $A_{\sigma}$ appearing in \eqref{HigherOrderNonlinearEstimate}
and \eqref{HigherOrderNonlinearDifferenceEstimate} is identical with
the function $A$ appearing in \eqref{NonlinearEstimate}
and \eqref{NonlinearDifferenceEstimate}. Recall that $A$ is assumed 
to be increasing, and that 
$\X^{\sigma}_{T} = \Lambda^{s-\sigma} \X^{s}_{T}$ by definition.

As in section \ref{Existence}, $\norm{u_{j}}_{\X^{s}_{T}} \le R$,
where $R = 2 C \bigl( \Sobnorm{f}{s} + \Sobnorm{g}{s-1} 
\bigr)$, $C$ is the constant appearing in \eqref{HomogeneousEstimate}
and $T > 0$ is chosen so small that $2 C_{T} A(R) \le 1$.
Another induction, using \eqref{HigherOrderNonlinearEstimate}
and \eqref{HomogeneousEstimate},  gives
$\norm{u_{j}}_{\X^{\sigma}_{T}} \le R'$, where
$R' = 2 C \bigl( \Sobnorm{f}{\sigma} + \Sobnorm{g}{\sigma-1} \bigr)$.
Furthermore, by \eqref{CauchyCondition} we have
$\norm{u_{j}- u_{j-1}}_{\X^{s}_{T}} \le 2^{-j} R$.

Thus \eqref{HigherOrderNonlinearDifferenceEstimate} 
yields $B_{j+1} \le B_{j}/2 + K2^{-j}$, where $B_{j} = \norm{u_{j}- 
u_{j-1}}_{\X^{\sigma}_{T}}$ and $K = C_{T} A(R') R$.
By induction, this implies that $B_{j} \le B_{0}2^{-j} + 2Kj2^{-j}$
for $j \ge 0$. Thus, $(u_{j})$ is Cauchy in $\X^{\sigma}_{T}$.
\subsection{Smooth dependence on initial data}\label{SmoothDependenceProof}
We prove the assertion in remark (3) following Theorem
\ref{AbstractWPTheorem}.

Fix $r \ge 1$. We prove $\delta \mapsto u_\delta \in \X^s_T$ is $C^r$.
By a finite time-step argument as in \ref{Persistence},
it suffices to prove this for a $T > 0$ which depends continuously on
$$
  E = \sup_{\delta \in I}
  \datanorm{f_\delta}{g_\delta}{s},
$$
where $I = [-\delta_0,\delta_0]$.

Denote the iterates by $u_j(\delta) \in \X^s_T$.
Thus, $u_{0}(\delta)$ solves the homogeneous wave equation
with initial data $(f_\delta,g_\delta)$, and for $j \ge 1$,
\begin{equation}\label{Iterates}
  u_{j+1}(\delta) = u_{0}(\delta) + S\bigl( u_j(\delta) \bigr),
\end{equation}
where $S = W \mathcal N$.
If we set $u_{-1} \equiv 0$, this is valid for $j \ge -1$.

Since $C^{r}\bigl(I,\X^s_T \bigr)$ is a Banach space
when equipped with the norm
$$
  \norm{u} = \sum_{0 \le k \le r} \sup_{\delta \in I}
  \bignorm{u^{(k)}(\delta)}_{\X^s_T},
$$
it suffices to show that $u_j$ is Cauchy in this norm
for some $T(E) > 0$.

By \eqref{HigherDerBd} and the mean value theorem,
\begin{equation}\label{HigherOrderDiff}
  \bignorm{S^{(k)}(u) - S^{(k)}(v)}_{(T)}
  \le B \bigl( \max \{ \norm{u}_{\X^{s}_{T}}, 
  \norm{v}_{\X^{s}_{T}} \} \bigr) \norm{u - v}_{\X^{s}_{T}}
\end{equation}
for $k = 1,\dots,r$, where $B = \max_{2 \le k \le r+1} B_k$.

Let $C$ be the constant in \eqref{HomogeneousEstimate},
choose $T$ so that $C_T A(2CE) = \half$, and set
$$
  E_k = \sup_{\delta \in I}
  \norm{(d/d\delta)^k
  (f_\delta,g_\delta)}_{(s)}
$$
for $k = 1,\dots,r$.

By induction, as in section \ref{Existence}, we have
\begin{equation}\label{ZeroOrderBd}
  \sup_\delta \norm{u_j(\delta)} \le 2CE.
\end{equation}
Taking one derivative in \eqref{Iterates}, we get
$$
  u_{j+1}' = u_0' + S'(u_j)(u_j'),
$$
so by \eqref{FirstDerBd},
$$
  \norm{u_{j+1}'(\delta)}_{\X^s_T}
  \le CE_1 + C_T A(2CE) \norm{u_j'(\delta)}_{\X^s_T},
$$
and since $C_T A(2CE) = \half$, we conclude that
$$
  \sup_\delta \norm{u_j'(\delta)} \le 2CE_1.
$$

Taking two derivatives in \eqref{Iterates} gives
$$
  u_{j+1}'' = u_0'' + S''(u_j)(u_j',u_j') + S'(u_j)(u_j'').
$$
Thus, using \eqref{HigherDerBd},
$$
  \norm{u_{j+1}''(\delta)}_{\X^s_T}
  \le CE_2 + B(2CE) (2CE_1)^2 + \half \norm{u_j''(\delta)}_{\X^s_T},
$$
so
$$
  \sup_\delta \norm{u_j''(\delta)} \le 2CE_2 + 2B(2CE) (2CE_1)^2.
$$
Continuing like this, one finds that for $k = 0,\dots,r$ and all $j$,
\begin{equation}\label{DerivativeBounds}
  \sup_\delta \bignorm{u_j^{(k)}(\delta)} \le C_k(E,E_1,\dots,E_k),
\end{equation}
where $C_k$ is some continuous function.

By \eqref{HigherOrderNonlinearDifferenceEstimate} and \eqref{ZeroOrderBd},
and the fact that $C_T A(2CE) = \half$, we have
$$
  \norm{u_{j+1}(\delta) - u_j(\delta)}_{\X^s_T}
  \le \half
  \norm{u_{j}(\delta) - u_{j-1}(\delta)}_{\X^s_T}
$$
for $j \ge 0$, and it follows by induction that
$$
  \sup_\delta \norm{u_j(\delta) - u_{j-1}(\delta)}_{\X^s_T}
  \le CE 2^{-j}.
$$

Next, since
$$
  u_{j+1}' - u_j' = S'(u_j)(u_j' - u_{j-1}')
  + \bigl[ S'(u_j) - S'(u_{j-1}) \bigr] (u_{j-1}')
$$
for $j \ge 0$,we have
\begin{multline*}
  \norm{u_{j+1}'(\delta) - u_j'(\delta)}_{\X^s_T}
  \le \half
  \norm{u_{j}'(\delta) - u_{j-1}'(\delta)}_{\X^s_T}
  \\
  + B(2CE) \norm{u_{j}(\delta) - u_{j-1}(\delta)}_{\X^s_T}
  C_1(E,E_1).
\end{multline*}
Thus, redefining $C_1$,
$$
  \norm{u_{j+1}'(\delta) - u_j'(\delta)}_{\X^s_T}
  \le \half
  \norm{u_{j}'(\delta) - u_{j-1}'(\delta)}_{\X^s_T}
  + C_1(E,E_1) 2^{-j},
$$
and it follows by induction\footnote{Here, and below, we use the
following induction argument. Suppose $B_j$ is a sequence
of non-negative numbers such that $B_{j+1} \le B_j/2 + P(j) 2^{-j}$
for $j \ge 0$, where $P$ is a polynomial with non-negative
coefficients. Then there is a polynomial $Q$ such that
$B_j \le Q(j) 2^{-j}$ for all $j$. In fact, one can take
$Q(t) = \int_0^t 2 P(r) \, dr + B_0$.} that
$$
  \sup_\delta \norm{u_j'(\delta) - u_{j-1}'(\delta)}_{\X^s_T}
  \le C_1(E,E_1) (1 + j) 2^{-j},
$$
where once more we have redefined $C_1$. 

Since
\begin{align*}
  u_{j+1}'' - u_j'' &= \bigl[ S''(u_j) - S''(u_{j-1}) \bigr] (u_j',u_j')
  \\
  &\quad
  + S''(u_{j-1}) (u_j' - u_{j-1}',u_j')
  + S''(u_{j-1}) (u_{j-1}',u_j' - u_{j-1}')
  \\
  &\quad
  + \bigl[ S'(u_j) - S'(u_{j-1}) \bigr] (u_j'')
  + S'(u_j)(u_j'' - u_{j-1}''),
\end{align*}
for $j \ge 0$, we have, redefining $C_2$,
\begin{multline*}
  \norm{u_{j+1}''(\delta) - u_j''(\delta)}_{\X^s_T}
  \le \half
  \norm{u_{j}''(\delta) - u_{j-1}''(\delta)}_{\X^s_T}
  \\
  + C_2(E,E_1,E_2)
  \left( \norm{u_{j}(\delta) - u_{j-1}(\delta)}_{\X^s_T}
  +  \norm{u_{j}'(\delta) - u_{j-1}'(\delta)}_{\X^s_T} \right).
\end{multline*}
Thus, redefining $C_2$ again,
$$
  \norm{u_{j+1}''(\delta) - u_j''(\delta)}_{\X^s_T}
  \le \half
  \norm{u_{j}''(\delta) - u_{j-1}''(\delta)}_{\X^s_T}
  + C_2(E,E_1,E_2) (1 + j) 2^{-j}.
$$
It follows by induction that
$$
  \sup_\delta \norm{u_{j}''(\delta) - u_{j-1}''(\delta)}_{\X^s_T}
  \le C_2(E,E_1,E_2) (1 + j^2) 2^{-j},
$$
where $C_2$ has been redefined yet again.

Continuing in this manner, one finds that for $k = 0,1,\dots,r$,
$$
  \sup_\delta
  \bignorm{u_{j}^{(k)}(\delta) - u_{j-1}^{(k)}(\delta)}_{\X^s_T}
  \le C_k(E,E_1,\dots,E_k) (1 + j^k) 2^{-j},
$$
where $C_k$ has been redefined.
\section{Proof of Theorem \ref{SpecialWPTheorem}}\label{ProofOfSpecialWPTheorem}
As noted already, (S1--5) of section \ref{AbstractLWP} are satisfied.
Since $\mathcal N$ is local in time and maps $\X^s$ into
$\Distribution(\R^{1+n})$, (N1) of Theorem
\ref{AbstractWPTheorem} holds; (N2--4) are
obviously satisfied.

Let $0 < T < 1$ and $u \in \X^s_T$. Then $u$ is an equivalence class in
$\X^s$ (section \ref{Main}), and we denote by $\widetilde u$ an arbitrary
representative in $\X^s$ of this equivalence class. By assumption,
$\mathcal N(\widetilde u)$ belongs to $\Y^{s,\varepsilon}$,
so by Theorem \ref{BasicConditionsTheorem}, there is a unique
$v \in \X^s_T$ which solves $\square v = \mathcal N(\widetilde u)$ on
$(0,T) \times \R^n$ with vanishing initial data at $t = 0$.
But $v$ is independent of the choice of representative $\widetilde u$ of $u$, by
uniqueness of solutions of \eqref{CauchyProblem} and the fact that $\mathcal N$
is local in time, so we may write $v = W \mathcal N(u)$.

Let us prove \eqref{NonlinearDifferenceEstimate}.
Let $u, v \in \X^s_T$, and let $\widetilde u, \widetilde v \in \X^s$
be any two representatives of $u$ and $v$ respectively.
By Theorem \ref{BasicConditionsTheorem} followed by
\eqref{EpsilonGenericDifferenceEstimate}, we have
\begin{align*}
  \norm{W \bigl( \mathcal N(u) - \mathcal N(v) \bigr)}_{\X^{s}_{T}}
  &\le C_{T,\varepsilon}
  \norm{\Lambda_{+}^{-1} \Lambda_{-}^{\varepsilon-1}
  \bigl( \mathcal N(\widetilde u) - \mathcal N(\widetilde v) \bigr) }_{\X^{s}}
  \\
  &\le  C_{T,\varepsilon}
  A\bigl( \max\{\norm{\widetilde u}_{\X^s},\norm{\widetilde v}_{\X^s} \} \bigr)
  \norm{\widetilde u- \widetilde v}_{\X^s}.
\end{align*}
We may assume that $A$ is increasing, whence
$$
  \norm{W \bigl( \mathcal N(u) - \mathcal N(v) \bigr)}_{\X^{s}_{T}}
  \le C_{T,\varepsilon}
  A\bigl( \norm{\widetilde w}_{\X^s} + \norm{\widetilde v}_{\X^s} \bigr)
  \norm{\widetilde w}_{\X^s},
$$
where $\widetilde w = \widetilde u - \widetilde v$.
Now pass to the limit as $\norm{\widetilde w}_{\X^s} \to \norm{u-v}_{\X^s_T}$
and $\norm{\widetilde v}_{\X^s} \to \norm{v}_{\X^s_T}$.
Thus \eqref{NonlinearDifferenceEstimate} holds with $A(R)$ replaced by
$A(3R)$.

Thus, the hypotheses of Theorem \ref{AbstractWPTheorem}
are satisfied, so (I--III) hold.

It remains to prove (IV) and (V). Let $W_{T} : \Y^{s,\varepsilon} \to \X^{s}$
be as in remark (2) following Theorem \ref{BasicConditionsTheorem}.
Given $(f,g)$, let $u_0$ be the solution of the homogeneous
wave equation with data $(f,g)$. As in the proof of Theorem
\ref{BasicConditionsTheorem},
\begin{equation}\label{GlobalHomogeneous}
  \norm{\chi(t) u_{0}}_{\X^{s}} \le C \datanorm{f}{g}{s}.
\end{equation}
Consider the sequence $u_j \in \X^s$ of iterates, given inductively by
$$
  u_{j} = \chi(t) u_{0} + W_{T} \mathcal N(u_{j-1})
$$

If $\mathcal N$ is $C^\infty$, then $S = W_T \mathcal N : \X^s \to \X^s$
is $C^\infty$, and the argument in section \ref{SmoothDependenceProof}
proves (IV).

Finally, we prove that \eqref{HigherOrderEstimate} implies property (V).
As shown in section \ref{Persistence}, it suffices to prove property 
(I$'$) stated therein. Fix $\sigma \ge s$. We may
assume that \eqref{HigherOrderEstimate}
holds with $A_{\sigma} = A$, where $A$ is the function appearing in
\eqref{EpsilonGenericDifferenceEstimate}.

Given $(f,g) \in H^{\sigma} \times H^{\sigma-1}$, choose $0 < T < 1$
so small that $2C_{T}A(R) \le 1$, where $R$ is twice the right hand 
side of \eqref{GlobalHomogeneous}.
A simple induction argument, essentially like the one
in section \ref{Existence}, reveals that $\norm{u_{j}}_{\X^{s}} \le R$
for all $j$ and that $(u_{j})$ is a Cauchy sequence in $\X^{s}$
whose limit $u$ solves \eqref{GenericSystem} on $(0,T) \times \R^{n}$
with initial data $(f,g)$.

Another induction, using \eqref{HigherOrderEstimate}, shows that
$\norm{u_{j}}_{\X^{\sigma}} \le R'$ for all $j$, where $R' = 2 C
\bigl( \Sobnorm{f}{\sigma} + \Sobnorm{g}{\sigma-1} \bigr)$.
Thus, by hypothesis (i) of Theorem \ref{BasicConditionsTheorem}, the
sequence of iterates is bounded in the Hilbert space $\scrH^{\sigma,\theta}$,
hence it converges weakly in that space, so it converges in the sense
of distributions to an element of $\scrH^{\sigma,\theta}$.
Thus, the limit $u \in \X^{s}$ must in fact belong to $\scrH^{\sigma,\theta}$,
whence $u \in C\bigl( [0,T],H^{\sigma} \bigr) \cap C^{1}\bigl( 
[0,T],H^{\sigma-1} \bigr)$.
\end{document}